%% file: main.tex
\documentclass[10pt]{article} 
\include{0Prae}

\title{Computation of the self-diffusion coefficient with low-rank tensor methods: application to the simulation of a cross-diffusion system.}
\usepackage{authblk}
\author[1,2]{Jad Dabaghi}
\author[1,2]{Virginie Ehrlacher}
\author[3]{Christoph Strössner}
\affil[1]{Ecole des Ponts ParisTech, Marne-la-Vall\'ee, France}
\affil[2]{INRIA Paris, France}
\affil[3]{Institute   of   Mathematics, EPF  Lausanne,  Switzerland}
\date{}
\begin{document}

\maketitle

{\centering \vspace{-0.8cm}
\url{jad.dabaghi@enpc.fr}, \url{virginie.ehrlacher@enpc.fr}, \url{christoph.stroessner@epfl.ch} \\ \vspace{1cm}
\centering \large{\date{April 4, 2022}} \\ \vspace{1cm}
}

\begin{abstract} 
Cross-diffusion systems arise as hydrodynamic limits of lattice multi-species interacting particle models. The objective of this work is to provide a numerical scheme for the simulation of the cross-diffusion system identified in [J. Quastel, \emph{ Comm. Pure Appl. Math.}, 45 (1992), pp. 623--679]. To simulate this system, it is necessary to provide an approximation of the so-called self-diffusion coefficient matrix of the tagged particle process. Classical algorithms for the computation of this matrix are based on the estimation of the long-time limit of the average mean square displacement of the particle. In this work, as an alternative, we propose a novel approach for computing the self-diffusion coefficient using deterministic low-rank approximation techniques, as the minimum of a high-dimensional optimization problem. \revision{The computed self-diffusion coefficient} is then used for the simulation of the cross-diffusion system using an implicit finite volume scheme. 
\end{abstract}

\input{1Introduction}
\input{2Setting}

\input{3AutoDiffCoefficient}

\input{4LowRankApproximation}

\input{5DeterministicSolving}

\input{6NumericalExperiments}

\input{7Conclusion}

\subsection*{Acknowledgements}
This work was initiated during the CEMRACS 2021 summer school and the authors would like to thank the CIRM for hosting them at this occasion. 
The authors would like to thank Mi-Song Dupuy, Daniel Kressner and Guillaume Enchéry for insightful discussions. 
The authors acknowledge support from the ANR project COMODO (ANR-19-CE46-0002), \revision{from the I-Site FUTURE} and from the Center on Energy and Climate Change (E4C).

\input{main.bbl}
\end{document}

%% file: 0Prae.tex
\usepackage{fullpage}
\usepackage[T1]{fontenc} 
\usepackage[latin1]{inputenc}
\usepackage[english]{babel}
\usepackage{selinput}
\SelectInputMappings{adieresis={ä},germandbls={ß}}

\usepackage{amsmath} 
\usepackage{amssymb} 
\usepackage{amsthm} 
\usepackage{enumerate}

\theoremstyle{definition}
\newtheorem{rmrk}{Remark}[section]

\usepackage{bm,bbm}
\usepackage{enumerate}
\usepackage{graphicx} 
\usepackage{epstopdf}
\usepackage[framed,numbered]{matlab-prettifier}
\usepackage{algorithm} 
\usepackage[noend]{algpseudocode}

\graphicspath{{./Figures/}}
\usepackage{multirow}
\usepackage{hyperref}
\usepackage{mathrsfs}
\usepackage{tikz}
\usepackage[normalem]{ulem}
\usepackage{cancel}
\usepackage{tabls}
\usetikzlibrary{shapes,arrows}

\newcommand{\R}{\mathbb{R}}
\newcommand{\bra}[1]{\left( #1 \right) }

\newcommand{\norm}[1]{\lVert #1 \rVert}
\newcommand{\abs}[1]{\left\lvert #1 \right\rvert}

\newcommand{\eps}{\varepsilon}

\DeclareMathOperator{\argmin}{argmin}

\newcommand{\umin}[1]{\underset{#1}{\min}\;}

\newcommand{\re}{\texttt{\textcolor{red}{red}}}
\newcommand{\bl}{\texttt{\textcolor{blue}{blue}}}

\newcommand{\matD}{\mathbb{D}}
\newcommand{\matA}{\mathbb{A}}
\newcommand{\matJ}{\mathbb{J}}

\newcommand{\matS}{\mathbb{S}}
\renewcommand{\bm}{\mathbf}
\newcommand{\nab}{\boldsymbol{\nabla}}
\newcommand{\bv}{\bm v}
\newcommand{\bx}{\bm x}
\newcommand{\bu}{\bm u}
\newcommand{\by}{\bm y}
\newcommand{\bz}{\bm z}
\newcommand{\bs}{\bm s}
\newcommand{\bw}{\bm w}
\newcommand{\bU}{\bm U}
\newcommand{\bn}{\bm n}
\newcommand{\bB}{\bm B}

\newcommand{\bEta}{\boldsymbol{\eta}}
\newcommand{\dps}{\displaystyle}

\newcommand{\Ip}{I_p}
\newcommand{\Nt}{P_t}
\newcommand{\dist}{\mathrm{dist}}
\newcommand{\Th}{\mathcal{T}_h}
\newcommand{\Tf}{T_f}

\newcommand{\Eh}{\mathcal{E}_h}
\newcommand{\Ehint}{\mathcal{E}_h^{\mathrm{int}}}
\newcommand{\Ehext}{\mathcal{E}_h^{\mathrm{ext}}}
\newcommand{\EK}{\mathcal{E}_K}

\newcommand{\hK}{h_K}

\newcommand{\Ne}{N_{\mathrm{e}}}

\newcommand{\F}{\mathcal{F}}

\newcommand{\epslin}{\varepsilon_{\mathrm{lin}}}

\newcommand{\revision}[1]{\textcolor{black}{#1}}

%% file: 1Introduction.tex
\section*{Introduction}

Cross-diffusion systems appear in various application fields such as population dynamics~\cite{Shigesada:Kaw:Teramato:1979}, tumor growth in medical biology~\cite{Jackson:Byrne:2002}, or diffusion processes in materials science. In particular, these models are used in order to simulate diffusion processes within mixture of chemical compounds, which occurs for instance during physical vapor deposition processes for the fabrication of thin film solar cells~\cite{Bakhta18,2010xlv}. Such cross-diffusion systems read as nonlinear systems of coupled degenerate parabolic partial differential equations describing the time evolution of diffusion processes within multi-component systems. The mathematical analysis of such cross-diffusion systems has recently attracted the interest of many mathematicians because it yields quite challenging new difficulties~\cite{Burger10,Jungel:2015,jungel2016entropy,bruna2012diffusion}.

Some cross-diffusion systems are derived as hydrodynamic limits of stochastic processes at the microscopic level, for instance of some lattice-based stochastic hopping models such as the ones studied in~\cite{Liggett85,Quastel92,Taylor16,Spohn91,Kipnis99,Blondel18,Redig20,Sasada10, komorowski2012fluctuations,Arita17,Arita18}.
In this work, we focus more specifically on cross-diffusion systems which read as the hydrodynamic limit of the multi-species symmetric exclusion process studied in~\cite{Quastel92,erignoux2016limite}.

One major difficulty of the numerical simulation of this system is that it requires the evaluation of the so-called self-diffusion coefficient matrix of the tagged particle process~\cite{Liggett85,Blondel18,Ertul21,Quastel92,Shapira20,erignoux2016limite}.

The most classical method to compute numerical approximations of these self-diffusion coefficients is to use a Monte Carlo scheme, since these coefficients can be expressed as the long-time limit of the time average expectation of the mean square deviation of the tagged particle~\cite{Quastel92}. 
However, this type of stochastic method is typically very slow to converge because of the high variance of the quantity the expectation of which has to be computed.

In this work, we propose a different approach exploiting the fact that the self-diffusion coefficients can be explicitly obtained using the unique solutions 
of infinite-dimensional deterministic minimization problems~\cite{Blondel18,Quastel92}. 
These problems can be approximated by high-dimensional finite-dimensional minimization problems, which however suffer from the curse of dimensionality. 
We mitigate this curse by using low-rank tensor methods in order to compute an approximate solution. This leads to a minimization problem \revision{over the set of low-rank tensors}, which we solve using a classical alternating scheme~\cite{Beylkin05,Grasedyck13,Rohwedder13}.
Our numerical experiments demonstrate that this low-rank approach leads to very accurate approximations of the self-diffusion coefficients. 
\revision{These results primarily serve as proof of concept and motivate the extension to more sophisticated low-rank approximation formats~\cite{Grasedyck13} in future work.}

Using this numerical approximation of the self-diffusion coefficients, the next step of this survey consists in computing the solution of the full cross-diffusion system. We then employ a particular cell-centered finite volume method, which satisfies local mass balance by construction~\cite{Eymard:Gallouet:Herbin:2000,Andreianov:Bendahmane:Ruiz2011,Shigesada:Kaw:Teramato:1979,Cances20}. 

This work is organized as follows. In Section~\ref{sec:setting:discrete}, we introduce the lattice-based stochastic hopping model and its hydrodynamic limit. Section~\ref{sec:auto:diff:coef} is dedicated to the presentation of the method we consider for the computation of the self-diffusion matrix from the resolution of some high-dimensional optimization problem. The low-rank approach mitigating the curse of dimensionality is presented in Section~\ref{sec:low:rank}.  In Section~\ref{sec:FiniteVolumes}, we introduce the cell-centered finite volume scheme for the resolution of \revision{a} cross-diffusion model \revision{inspired by~\cite{Quastel92}}. Finally, in Section~\ref{sec:numerical:experiments}, we perform numerical simulations for this problem.

%% file: 2Setting.tex
\section{Hydrodynamic limit of a lattice-based stochastic hopping model}\label{sec:setting:discrete} 

\begin{figure}
    \centering
    \includegraphics[width=0.9\textwidth]{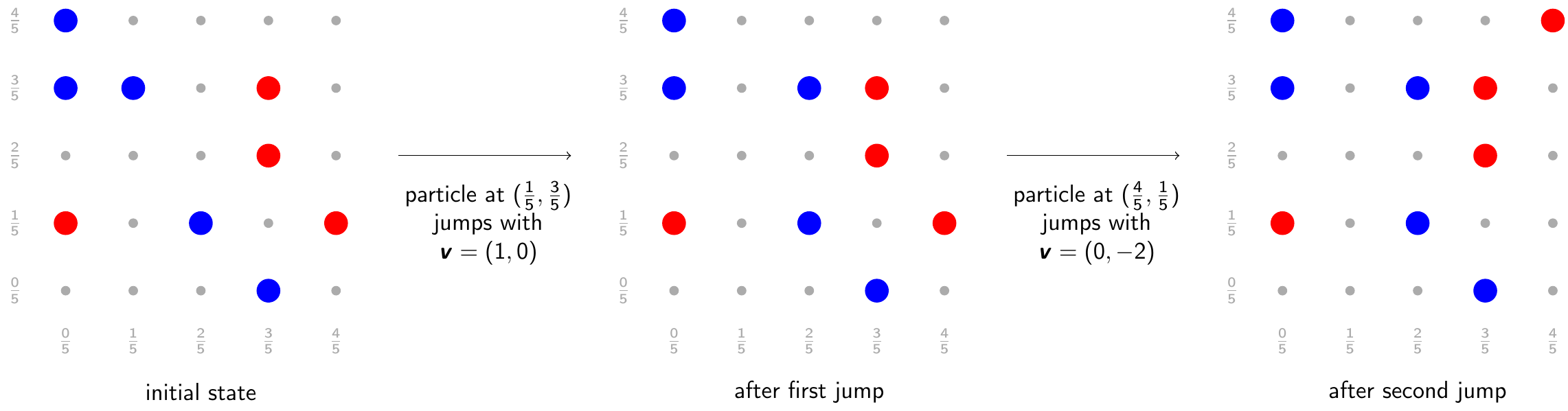}
\caption{The evolution of state of the particles on a lattice with $n=5,d=2$. Left: Initial state. Middle: After the blue particle at $\bs = (\frac{1}{5},\frac{3}{5})$ jumped with $\mathbf v=(1,0)$ to $\bs + \frac{1}{5} \mathbf v$. Right: After the red particle at position $(\frac{4}{5},\frac{1}{5})$ jumped with $\mathbf v=(0,-2)$. Since the lattice is periodic the particle ends up in position $(\frac{4}{5},\frac{4}{5})$.}
\label{fig:Hopping}
\end{figure}

\revision{Our motivation for this work stems from} the {lattice-based stochastic hopping model} proposed \revision{in}~\cite{Quastel92}, which describes the evolution of a mixture of multi-species particles, the positions of which are clamped on a (periodic) lattice {$T_n:=\{0,\dots,\frac{n-1}{n}\}^d$}, $d=1,2,3$. 
After a certain amount of time, a single particle jumps to a nearby node according to a given Markovian random walk.  
Let $\bv_k{:= ( (\bv_k)_i )_{1\leq i \leq d}} \in {\mathbb{Z}}^d\setminus\{\bm 0\}$
for $k = 1,\dots,K$ be the displacement vectors for all possible jumps that {one} particle can make, i.e. a particle at position $\bx \in T_n$
lands at position $\bx+{\frac{1}{n}}\bv_k \in {T_n}$. 
The {rate of } jumping in the direction $\frac{1}{n}\bv_k$ is given by $p_k \in \left]0,1\right]$. We {assume that } the particles can be of two possible species, either species $\bl$ or species $\re$. We study how the distribution of the { particles belonging to each} species {evolves} over time in the limit when the number of lattice points and  particles tends to infinity. See Figure~\ref{fig:Hopping} for a {schematic illustration of this lattice-hopping process}.

\medskip

In {this work}, we are interested in simulating the hydrodynamic limit of this lattice-based hopping model, which was identified in~\cite{Quastel92,Blondel18}. 
In the limit {$n\to +\infty$}, {this hydrodynamic limit reads as a cross-diffusion system defined on} the $d$-dimensional torus $T^d{:= (\mathbb{R}/\mathbb{Z})^d}$. 
\revision{At this scale, we consider} densities $\rho_{\re}: T^d \times \mathbb{R}_{+} \to [0,1]$ and $\rho_{\bl}: T^d \times \mathbb{R}_{+}  \to [0,1]$\revision{, where} $\rho_{\re}(\bx,t)$ (respectively  $\rho_{\bl}(\bx,t)$) denotes the local volumic fraction of $\re$ (respectively $\bl$) particles at point $\bx \in T^d$ and time $t\geq 0$. 
The local volumic fractions $\rho_{\re}$ and $\rho_{\bl}$ {are shown to be} solutions to the following cross-diffusion system~\cite{Quastel92}
\begin{align}
\label{eq:cross:diff:system}
    \frac{\partial}{\partial t} 
    \begin{pmatrix} 
    \rho_{\re} \\ \rho_{\bl} 
    \end{pmatrix} 
    = \frac{1}{2} \nabla \cdot \left[
    \begin{pmatrix} 
    \dps \frac{\rho_{\bl}}{\rho} \matD_s(\rho) + \dps \frac{\rho_{\re}}{\rho} \matD & \dps \frac{\rho_{\re}}{\rho} (\matD-\matD_s(\rho)) \\ 
    \dps \frac{\rho_{\bl}}{\rho} (\matD-\matD_s(\rho)) 
    & \dps \frac{\rho_{\re}}{\rho} \matD_s(\rho) + \dps \frac{\rho_{\bl}}{\rho} \matD 
    \end{pmatrix} 
    \begin{pmatrix} \nabla \rho_{\re} \\ \nabla \rho_{\bl}
    \end{pmatrix}\right], \quad {\mbox{ in }T^d \times \mathbb{R}_+^*}, 
\end{align}
{where $\rho:= \rho_\re + \rho_\bl$.} 
We complement system~\eqref{eq:cross:diff:system} by the initial conditions $\rho_{\re}(0,\bx) = \rho_{\re}^0(\bx),\ \rho_{\bl}(0,\bx) = \rho_{\bl}^0(\bx)$, where 
{$\rho_{\re}^0, \rho_\bl^0 \in L^\infty(T^d)$ satisfying} $0 \leq \rho_{\re}^0+\rho_{\bl}^0 \leq 1$ {almost everywhere in $T^d$}. 
{In this case, it can be shown~\cite{Quastel92} that there exists a unique solution $(\rho_\re, \rho_\bl) \in L^\infty(T^d \times \mathbb{R}_+)^2$ to system \ref{eq:cross:diff:system}. In addition, it holds that for almost all $t\geq 0$ and $\bx\in T^d$, $0 \leq \rho_\re(\bx,t)$, $0\leq \rho_\bl(\bx, t)$ and 
$\rho_\re(\bx,t) + \rho_\bl(\bx,t) \leq 1$.}

The matrix $\matD {:= (\matD_{ij})_{1\leq i,j \leq d}} \in \R^{d \times d}$ is defined so that $\matD_{ij} := \sum_{k=1}^K p_k (\bv_k)_i (\bv_k)_j$  for all $1 \leq i,j \leq d$. 
The self-diffusion matrix application
$$
\matD_s : \left\{ 
\begin{array}{ccc}
[0,1] & \to & \mathbb{R}^{d\times d}\\
\overline{\rho} & \mapsto & \matD_s(\overline{\rho}):= \left(\matD_{s,ij}(\overline{\rho})\right)_{1\leq i,j \leq d}\\
\end{array}\right.
$$
will be defined in more details in Section~\ref{sec:auto:diff:coef}, (see also~\cite{Blondel18,Ertul21}) and is such that for all $\overline{\rho}\in [0,1]$, $\matD_s(\overline{\rho})$ is a symmetric positive semi\revision{-}definite matrix. 
The system~\eqref{eq:cross:diff:system} is a cross-diffusion system and is a highly nonlinear degenerate parabolic system. \revision{Numerical simulations of system~(\ref{eq:cross:diff:system}) require as a first step the computation of numerical approximations of $\matD_s(\overline{\rho})$ for any value $\overline{\rho}\in [0,1]$, and we propose in this work a novel numerical approach for this task.}

%% file: 3AutoDiffCoefficient.tex
\section{Self-diffusion matrix} \label{sec:auto:diff:coef}

\subsection{Definition of the self-diffusion matrix}

For all $\bu \in \R^d$ and $\overline{\rho}\in [0,1]$, the self-diffusion coefficient $\bu^T\mathbb{D}_s\left(\overline{\rho}\right)\bu$ is obtained by considering the so-called tagged particle process on 
$S:= \mathbb{Z}^d \setminus \{\bm 0\}$. Let us first introduce some notation. Let $\bEta:= (\eta_\bs)_{\bs\in S} \in \{0,1\}^S$. For all $\by\neq \bz \in S$, we define by $\bEta^{\by, \bz}:= (\eta^{\by, \bz}_\bs)_{\bs\in S} \in \{0,1\}^S$ so that
$$
\eta^{\by, \bz}_\bs:=\left\{
\begin{array}{ll}
\eta_\bs & \mbox{if } \bs\neq \by, \bz,\\
\eta_\by& \mbox{ if }\bs = \bz,\\
\eta_\bz & \mbox{ if } \bs = \by,\\
\end{array}\right. 
$$
and for all $\bw \in S$, we define by $\bEta^{\bm 0, \bw}:= (\eta^{\bm 0, \bw}_\bs)_{\bs\in S} \in \{0,1\}^S$ so that
$$
\eta^{\bm 0, \bw}_\bs:=\left\{
\begin{array}{ll}
\eta_{\bs+\bw} & \mbox{if } \bs\neq -\bw, \\
0& \mbox{ if }\bs = -\bw.\\
\end{array}\right. 
$$

More precisely, denoting by $H:= \left\{ \Psi: \{0,1\}^S \to \mathbb{R}\right\}$, it then holds that
\begin{equation}\label{eq:infinite}
\bu^T \matD_s(\overline{\rho}) \bu:= \revision{2} \mathop{\inf}_{\Psi\in H} \mathbb{E}_{\overline{\rho}^\otimes}\left[\sum_{k=1}^K p_k
\bra{ (1- \eta_{\bv_k}) \bra{
\bu\cdot \bv_k + \Psi(\bEta^{\bm 0,\bv_k}) - \Psi(\bEta)}^2 + \revision{\frac{1}{2}}
\sum_{\substack{\by \in S \setminus \{ {\bm 0} \} \\ \by+\bv_k \neq \bm 0 } } \bra{
\Psi(\bEta^{\by+\bv_k,\by})-\Psi(\bEta)}^2 }  \right], 
\end{equation}
where the notation $\mathbb{E}_{\overline{\rho}^\otimes}$ refers to the fact that the expectation is computed on all random variables $\bEta:= (\eta_\bs)_{\bs\in S}$ so that the random variables $\eta_\bs$ are \revision{independently} identically distributed random variables according to a Bernoulli law with parameter $\overline{\rho}$. Problem (\ref{eq:infinite}) thus reads as an infinite-dimensional optimization problem which we are going to discretize as follows. 

\medskip
\begin{rmrk} 
Naturally, for all $\overline{\rho}\in [0,1]$, since $\matD_s(\overline{\rho})$ is a symmetric matrix, one can easily know the full matrix $\matD_s(\overline{\rho})$ from the knowledge of $\bu^T \matD_s(\overline{\rho})\bu$ for a few vectors $\bu \in \mathbb{R}^d$. 
\end{rmrk}

\subsection{Finite-dimensional approximation of the self-diffusion matrix}

\begin{figure}
    \centering
    \includegraphics[width=0.9\textwidth]{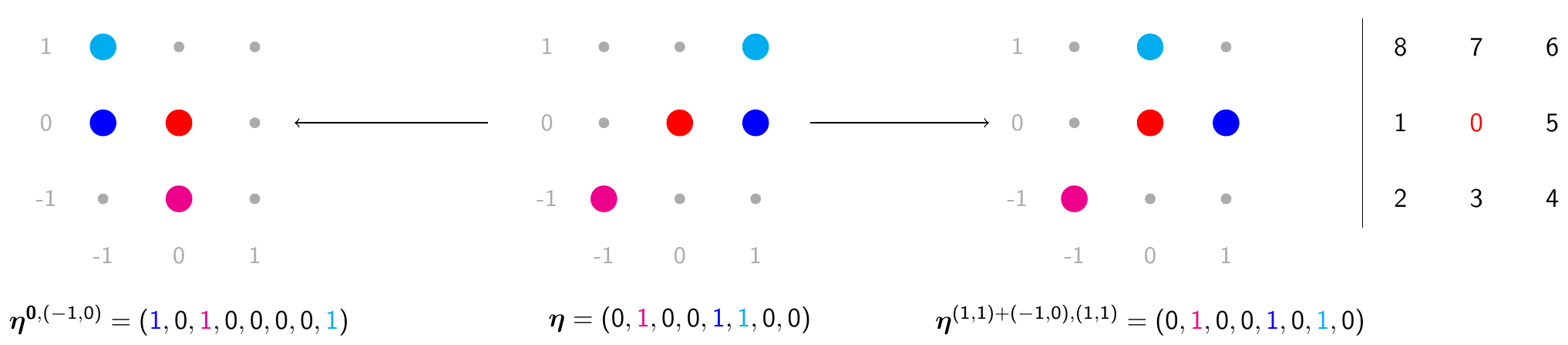}
    \caption{
    Middle: Visualization of one particular $\bEta \in S_1$ for $d=2$ with three occupied sites marked in blue, cyan and magenta. Additionally, we mark $\bm 0$ in red.
    Left: Visualization of the occupied sites of $\bEta^{\bm 0, (-1,0)}$. This can be seen as a jump of the imaginary red particle one step to the left, followed by immediately relabeling of the sites such that red particle remains at $\bm 0$. By exploiting the periodicity, we obtain $\bEta^{\bm 0, (-1,0)} \in S_1$.
    Right: Visualization of the occupied sites of $\bEta^{(1,1)+(-1,0),(1,1)}$. This can be seen as jump of the cyan particle one step to the left.
    Far right: The enumeration scheme for the sites. 
    }
    \label{fig:Jumps}
\end{figure}

Let $M \in \mathbb{N}^*$ denote a discretization parameter and introduce the finite grid $S_M:= \{-M, \cdots, M\}^d \setminus \{\bm 0\}$.  For any $\bEta \in \{0,1\}^{S_M}$, we can construct by periodicity an extension $\widetilde{\bEta}:= (\widetilde{\eta}_\bs)_{\bs\in S} \in \{0,1\}^S$ by assuming with a slight abuse of notation that the site $\bm 0$ is occupied, i.e. $\widetilde{\eta}_\bs = 1$ for $\bs \in (2M+1)\cdot\mathbb{Z}^d \setminus \{\bm 0\}$. Using this notation, for all $\bEta\in \{0,1\}^{S_M}$ and all $\by, \bz, \bw \in S$, we define $\bEta^{\by, \bz}\in \{0,1\}^{S_M}$ and $\bEta^{\bm 0, \bw}\in \{0,1\}^{S_M}$ as 
$$
\bEta^{\by, \bz}:= \left( \widetilde{\eta}^{\by, \bz}_\bs \right)_{\bs \in S_M} \quad \mbox{ and }\quad \bEta^{\bm 0, \bw}:= \left( \widetilde{\eta}^{\bm 0, \bw}_\bs \right)_{\bs \in S_M}.
$$
Let us then denote by $N:= (2M+1)^d -1 = {\rm Card}(S_M)$. Then, any $\bEta = (\eta_\bs)_{\bs\in S_M} \in S_M$ can be equivalently viewed as an element $\bEta := (\eta_i)_{1\leq i \leq N}\in \{0,1\}^N$ by enumerating the different sites of $S_M$. Figure~\ref{fig:Jumps} shows an illustration of $\bEta^{\by, \bz}$ and $\bEta^{\bm 0,\bw}$.

For all $\ell \in \{0,\dots,N\}$, let $C_{M,\ell} := \{ \bEta \in \{0,1\}^{S_M} | \sum_{\bs \in S_M} \eta_\bs = \ell\}$ denote the set of all possible configurations of the particles on $S_M$ so that the total number of occupied sites is equal to $\ell$. 


Let us also denote by $H_M:= \left\{ \Psi: \{0,1\}^{S_M} \to \mathbb{R}\right\}$. For every $\bu \in \R^d$ and $\revision{\ell} \in \{0,\dots,N\}$, we then introduce  the quadratic functional $A_{M,\ell}^\bu : H_M \to \mathbb{R}$ defined by
\begin{equation*}
A_{M,\ell}^{\bm u}(\Psi) := \frac{1}{|C_{M,\ell}|} \sum_{\bEta \in C_{M,\ell}} 
\sum_{k=1}^K p_k
\bra{ (1- \eta_{\bv_k}) \bra{
\bu\cdot \bv_k + \Psi(\bEta^{\bm 0,\bv_k}) - \Psi(\bEta)}^2 + \revision{\frac{1}{2}}
\sum_{\substack{\by \in S_M \\ \by+\bv_k \neq \bm 0 } } \bra{
\Psi(\bEta^{\by+\bv_k,\by})-\Psi(\bEta)}^2   }.
\end{equation*}

Then, assuming that $\overline{\rho} = \frac{\ell}{N}$ for some $0\leq \ell \leq N$, one can define~\cite{Blondel18} for all $\bu \in \mathbb{R}^d$, 
\begin{equation}
\label{eq:approx:Ds}
\bu^T \mathbb{D}^M_s\left(\frac{\ell}{N}\right) \bu:= \revision{2} \umin{\Psi \in H_M} A_{M, \ell}^{\bm u} (\Psi).
\end{equation}
\revision{Note that (\ref{eq:approx:Ds}) can then be seen as a finite-dimensional approximation of optimization problem (\ref{eq:infinite}). Indeed, it }is proved in~\cite{landim2002finite} that $\displaystyle \mathop{\lim}_{\begin{array}{c}
    M\to +\infty \\
    \frac{\ell}{N} \to \overline{\rho} \\
    \end{array}}
\bu^T \mathbb{D}^M_s\left(\frac{\ell}{N}\right) \bu = \bu^T \mathbb{D}_s(\overline{\rho}) \bu$. 

\subsection{\revision{Combined} minimization problem}\label{sec:relaxing}
The collection of sets $C_{M,0},\dots,C_{M,N}$ form a partition of the set $\{0,1\}^{S_M}$. Observe for a given $\Psi\in H_M$,$A_{M,\ell}^{\bm u}(\Psi)$ only depends on the values of $\Psi(\bEta)$ for $\bEta \in C_{M,\ell}$, since $\bEta \in C_{M,\ell}$ implies that also $\bEta^{\bm 0,\bm v_k} \in C_{M,\ell}$ and $\bEta^{\bm y + \bm v_k, \bm y} \in C_{M,\ell}$ for all $1\leq k \leq K$.
As a consequence,  if $\Psi^{M,\bm u}_{\textsf{opt}}\in H_M$ is a minimizer of 
\begin{equation}\label{eq:minopt}
\min_{\Psi \in H_M} A^{\bm u}_M (\Psi),
\end{equation}
where 
\begin{equation}
    \label{eq:def:Au}
A_M^{\bm u}(\Psi) := \sum_{\bEta \in \{0,1\}^{S_M}} 
\sum_{k=1}^K p_k
\bra{ (1- \eta_{\bv_k}) \bra{
\bu\cdot \bv_k + \Psi(\bEta^{\bm 0,\bv_k}) - \Psi(\bEta)}^2 + \revision{\frac{1}{2}}
\sum_{\substack{\by \in S_M \\ \by+\bv_k \neq \bm 0 } } \bra{
\Psi(\bEta^{\by+\bv_k,\by})-\Psi(\bEta)}^2 },
\end{equation}
then it holds that $A_{M,\ell}^{\bm u}(\Psi^{M,\bm u}_{\textsf{opt}}) = \min_{\Psi \in H_M} A_{M,\ell}^{\bm u}(\Psi)$ for all $\ell \in \{0,\dots,N\}$. The knowledge of $\Psi^{M,\bm u}_{\textsf{opt}}$ then allows us to compute $\bu^T \mathbb{D}_s^M\left(\frac{\ell}{N}\right) \bu$ for all $0\leq \ell \leq N$ as $\revision{2} A_{M,\ell}^{\bm u} (\Psi^{M,\bm u}_{\textsf{opt}})$. Note that the minimization problem (\ref{eq:minopt}) is then independent of $\ell$, in contrast to (\ref{eq:approx:Ds}). 

\medskip 

For small values of $M$, one can compute  $\min_{\Psi \in H_M} A_M^{\bm u} (\Psi)$ explicitly. 
Indeed, the function $\Psi \in H_M$ can be equivalently identified as a tensor $\bm \Psi:= \left( \bm \Psi_{\bEta} \right)_{\bEta \in\{0,1\}^{S_M} } \in \R^{2^N}$,
where for all $\bEta\in \{0,1\}^{S_M}$, $\bm \Psi_{\bEta} = \Psi(\bEta)$. 
We see the quadratic functional $A_M^\bu$ as sum of $Q = NK \cdot 2^{N-1}$ quadratic terms, which can be expressed as $\norm{\mathbb{M}_M \bm \Psi + \bm b_M^\bu}_2^2$ for some matrix $\mathbb{M}_M \in \mathbb{R}^{2^N \times Q}$ and vector $\bm b_M^\bu \in \mathbb{R}^{2^N}$. Thus, (\ref{eq:minopt}) then boils down to solving the least squares problem
 \begin{equation}
    \label{eq:linlsq}
\umin{\bm \Psi \in \R^{2^N}} \norm{\mathbb{M}_M \bm \Psi + \bm b_M^\bu}_2^2.
\end{equation}
In the case where $d=2$, $M=1$ and $K=4$, it holds that $N=8$ so that $2^N = 256$ and $Q = 4096$. A solution of~\eqref{eq:linlsq} is then computed in practice up to a very high precision using $\texttt{lsqr}$ in Matlab. 

\medskip

However, for larger values of $M$ (i.e. larger values of $N$) this approach quickly becomes intractable, since the size of the matrix grows exponentially in $N$. The goal of Section~\ref{sec:low:rank} is to propose a low-rank approximation method in order to compute a numerical approximation of $\Psi^{M,\bm u}_{\textsf{opt}}$, hopefully in cases with $M$ large, which allows us to approximate $\bu^T \mathbb{D}_s^M\left(\frac{\ell}{N}\right) \bu$ for all $0\leq \ell \leq N$ by evaluating $A_{M,\ell}^{\bm u} (\Psi^{M,\bm u}_{\textsf{opt}})$.

\begin{rmrk}
Note that we could compute the optimal value of $\min_{\Psi \in H_M} A_{M,\ell}^{\bm u} (\Psi)$ directly, by solving a linear least squares problem with $\binom{N}{\ell}$ degrees of freedom (see Section~\ref{sec:relaxing}). While this is cheaper than solving the system with $2^N$ degrees of freedom, it still is intractable for many choices of $l$ for larger $N$.
\end{rmrk}

%% file: 4LowRankApproximation.tex
\section{Separable low-rank approximation}\label{sec:low:rank}

A function $R \in H_M$ is called {a separable or pure tensor product function} when it can be written as \[
R (\bEta) = \Pi_{\bs\in S_M} R_\bs(\eta_\bs), \quad {\forall \bEta = (\eta_\bs)_{\bs\in S_M}\in \{0,1\}^{S_M},}
\]
for some $R_\bs: \{0,1\} \to \R$ for $\bs\in S_M$. {Let $\mathcal{T}_M^1 \subset H_M$ denote the set of pure tensor product functions of $H_M$. }
It holds that 
\begin{equation*}
    \umin{\Psi \in H_M} A^{\bm u}_M(\Psi) \leq \umin{R \in \mathcal T^1_M} A^{\bm u}_M(R) .
\end{equation*}

The aim of this section is to present the numerical strategy we developed in this work so as to compute a numerical approximation \revision{$R^{\bu}_{\rm opt}$ of a minimizer of 
\begin{equation}\label{eq:approxmin}
    \mathop{\min}_{R \in \mathcal T_M^1} A^{\bm u}_M(R).
\end{equation}
Note that there always exists at least one minimizer to (\ref{eq:approxmin}) but that uniqueness is not guaranteed in general}. \revision{The approximation $R^{\bu}_{\rm opt}$} is then used in turn as an approximation of $\Psi^{\bm u}_{\textsf{opt}}$. We will present numerical tests in Section~\ref{sec:low:rank:error:experiments} which demonstrate that \revision{$R^{\bu}_{\rm opt}$ is indeed close to $\Psi^{\bm u}_{\textsf{opt}}$.}

\begin{rmrk}
\revision{
The goal of this work, is to demonstrate that optimization over the space of pure tensor product functions yields good minimizers of $A_M^{\bm u}$. The extension of this approach to to sums of pure tensor product functions or other low-rank approximation formats~\cite{Grasedyck13} offers the potential to yield even better minimizers and is subject to future work.}
\end{rmrk}

\subsection{Exploiting separability}\label{sec:rewritingAu}
Let us explain in this section why it is advantageous to minimize $A_M^{\bu}$ over the set $\mathcal T_M^1$. Indeed, let us rewrite $A^{\bu}_M(\Psi)$ for all $\Psi \in H_M$ as
\begin{align*}
    A^{\bm u}_M (\Psi) = \sum_{k=1}^K p_k \Bigg(& \sum_{\bEta \in \{0,1\}^{S_M}} (1-\eta_{\bv_k}) \Psi(\bEta^{\bm 0,\bv_k})^2 + \sum_{\bEta \in \{0,1\}^{S_M}} (1-\eta_{\bv_k}) \Psi(\bEta)^2 -2 \sum_{\bEta \in \{0,1\}^{S_M}} (1-\eta_{\bv_k}) \Psi(\bEta^{\bm 0,\bv_k})\Psi(\bEta) \nonumber \\ &+ 2(\bm u\cdot \bm  v_k) \bigg( \sum_{\bEta \in \{0,1\}^{S_M}} (1-\eta_{\bv_k}) \Psi(\eta^{\bm 0,\bv_k}) - \sum_{\bEta \in \{0,1\}^{S_M}} (1-\eta_{\bv_k}) \Psi(\bEta) \bigg) + 2^{N-1}(\bm u \cdot \bm v_k)^2  \\ &+ \revision{\frac{1}{2}}
    \sum_{\substack{\bm y \in S_M \\ \bm y+\bm v_k \neq \bm 0}} \bigg( \sum_{\bEta \in \{0,1\}^{S_M}} \Psi(\bEta^{\bm y+\bm v_k,\bm y})^2 + \sum_{\bEta \in \{0,1\}^{S_M}} \Psi(\bEta)^2 -
    2 \sum_{\bEta \in \{0,1\}^{S_M}} \Psi(\bEta)\Psi(\bEta^{\bm y+\bm v_k,\bm y}) \bigg) \nonumber
    \Bigg) .
\end{align*}
Note that for all $1\leq k \leq K$, there exists a bijection $\pi^{\bm 0, \bv_k}: S_M \to S_M$ such that for all $\bEta = (\eta_\bs)_{\bs\in S_M}$, 
$$
\bEta^{\bm 0, \bv_k} = \left( \eta_{\pi^{\bm 0, \bv_k}(\bs)}\right)_{\bs\in S_M}. 
$$
Similarly, for all $1\leq k \leq K$ and all $\by \in S_M$ such that $\by + \bv_k \neq \bm 0$, there exists a bijection $\pi^{\by, \by + \bv_k}: S_M \to S_M$ such that for all $\bEta = (\eta_\bs)_{\bs\in S_M}$, 
$$
\bEta^{\by, \by + \bv_k} = \left( \eta_{\pi^{\by, \by + \bv_k}(\bs)}\right)_{\bs\in S_M}. $$
Let us denote by $\sigma^{\bm 0, \bv_k}$ (respectively $\sigma^{\by, \by +\bv_k}$) the inverse of $\pi^{\bm 0, \bv_k}$ (respectively $\pi^{\by, \by+\bv_k}$).

Then, it holds that for all separable function $R \in \mathcal T_M^1$, so that $R(\bEta) = \Pi_{\bs\in S_M} R_\bs(\eta_\bs)$ for all $\bEta = (\eta_\bs)_{\bs\in S_M}$ for some $R_\bs:\{0,1\} \to \mathbb{R}$, 
\begin{align*}
  \displaystyle  A^{\bm u}_M (R) &  \displaystyle  = \sum_{k=1}^K p_k \Bigg[  \left( \Pi_{\bs\in S_M\setminus\{\bv_k\}} \Sigma_{\eta_\bs \in \{0,1\}} R_{\sigma^{\bm 0, \bv_k}(\bs)}(\eta_\bs)^2 \right) R_{\sigma^{\bm 0, \bv_k}(\bv_k)}(0)^2   + \left(\Pi_{\bs\in S_M\setminus\{\bv_k\}} \Sigma_{\eta_\bs \in \{0,1\}} R_{\bs}(\eta_\bs)^2 \right) R_{\bv_k}(0)^2 \\
  &+ 2(\bm u\cdot \bm  v_k) \bigg( \left(\Pi_{\bs\in S_M\setminus\{\bv_k\}} \Sigma_{\eta_\bs \in \{0,1\}} R_{\sigma^{\bm 0, \bv_k}(\bs)}(\eta_\bs)\right) R_{\sigma^{\bm 0, \bv_k}(\bv_k)}(0)  - \left(\Pi_{\bs\in S_M\setminus\{\bv_k\}} \Sigma_{\eta_\bs \in \{0,1\}} R_{\bs}(\eta_\bs) \right) R_{\bv_k}(0) \bigg)\\
  &+ \revision{\frac{1}{2}} \sum_{\substack{\bm y \in S_M \\ \bm y+\bm v_k \neq \bm 0}} \bigg[ \left(\Pi_{\bs\in S_M} \Sigma_{\eta_\bs \in \{0,1\}} R_{\sigma^{\by, \by +\bv_k}(\bs)}(\eta_\bs)^2 \right) + \left(\Pi_{\bs\in S_M} \Sigma_{\eta_\bs \in \{0,1\}} R_{\bs}(\eta_\bs)^2 \right) \\ & -2 \left(\Pi_{\bs\in S_M} \Sigma_{\eta_\bs \in \{0,1\}} R_{\bs}(\eta_\bs)R_{\sigma^{\by, \by+\bv_k}(\bs)}(\eta_\bs) \right) \bigg] + 2^{N-1}(\bm u \cdot \bm v_k)^2
    \Bigg] .\\
\end{align*}
Note that all these terms can be evaluated in $\mathcal{O}(N)$ operations, which avoids the summation over all $2^N$ possible $\bEta \in \{0,1\}^{S_M}$. Thus, the separability property allows us to efficiently evaluate $A^{\bm u}_M(R)$ for all $R\in \mathcal T_M^1$.

\subsection{Alternating least squares}

In this work, we use the classical Alternating Least Squares algorithm~\cite{Carroll70,Kolda09, de2008decompositions, oseledets2018alternating} in order to compute an approximate solution of the minimization problem $\min_{R \in \mathcal T_M^1} A_M^{\bm u}(R)$. 

The main idea is to find an approximation of $R^{\bm u}_{\textsf{opt}}$ by an iterative scheme which amounts to solving a sequence of small-dimensional linear problems. We start from an initial $R(\bEta) := \Pi_{\bs\in S_M}R_\bs(\eta_\bs)$. The first least squares problem is obtained by minimizing $A^{\bm u}(R)$ only with respect to a selected $R_{\bs_0}:\{0,1\} \to \R$ for some $\bs_0\in S_M$ leaving the other $R_\bs$, $\bs\neq \bs_0$ fixed. By partially evaluating  $A^{\bm u}(R)$ for all terms not depending on $R_{\bs_0}$, we obtain that $\min_{R_{\bs_0} \in \{\{0,1\} \to \R\}}A^{\bm u} (R)$ with $R(\bEta):= R_{\bs_0}(\eta_{\bs_0}) \Pi_{\bs\in S_M\setminus \{\bs_0\}} R_\bs(\eta_\bs)$ is equivalent to a quadratic optimization problem
\begin{equation*}
   \min_{R_{\bs_0} \in \{\{0,1\} \to \R\}} \alpha_1 R_{\bs_0}(1)^2 + \alpha_2 R_{\bs_0}(0)^2 + \alpha_3 R_{\bs_0}(1)R_{\bs_0}(0) + \alpha_4 R_{\bs_0}(1) + \alpha_5 R_{\bs_0}(0) + \alpha_6,
\end{equation*}
with constants $\alpha_1,\dots,\alpha_6 \in \R$ depending on the fixed $R_\bs$, $\bs\neq \bs_0$. This quadratic optimization problem always admits a unique optimal $R_{s_0}$, which is given by $R_{\bs_0}(1) = a$ and $R_{\bs_0}(0) = b$, where $a,b\in R$ are the solution of the linear system \[
\begin{pmatrix}
2\alpha_1 & \alpha_3 \\ \alpha_3 & 2 \alpha_2
\end{pmatrix}
\begin{pmatrix}
a \\ b
\end{pmatrix} =
\begin{pmatrix}
-\alpha_4 \\ - \alpha_5
\end{pmatrix}.
\]

This allows us to optimize $A^{\bm u}_M(R)$ with respect to individual $R_{\bs_0}$. 
By alternating the selected $\bs_0\in S_M$, we obtain the alternating least squares algorithm, which is formalized in Algorithm~\ref{alg:ALS}. 

\begin{algorithm}
\caption{Alternating least squares}\label{alg:ALS}
\begin{algorithmic}[1]
\algdef{SE}[SUBALG]{Indent}{EndIndent}{}{\algorithmicend\ }%
\algtext*{Indent}
\algtext*{EndIndent}
\State \textbf{Input}: initial functions $R^0_\bs: \{0,1\} \to \mathbb{R}$ for $\bs\in S_M$, vector $\bm u$, tolerance $\eps$
\State \textbf{Output}: approximation ${R}_{\rm opt}(\bEta) = \Pi_{\bs\in S_M}R_\bs(\eta_\bs)$ of $\argmin_{R \in \mathcal T_M^1} A^{\bm u}_M(R)$
\State $v_{\textsf{old}} = \infty$, $v_{\textsf{new}} = A_M^{\bm u}(R^0)$ with $R^0(\bEta):= \Pi_{\bs\in S_M} R^0_\bs(\eta_\bs)$
\State $\forall \bs\in S_M, \; R_\bs := R_\bs^0$. 
\While $\abs{v_{\textsf{old}}-v_{\textsf{new}}} > \eps \abs{v_{\textsf{new}}}$.
\State $v_{\textsf{old}} = v_{\textsf{new}}$
\For $\bs_0\in S_M$
\State $R_{\bs_0} = \argmin_{\widetilde{R}_{\bs_0} :\{0,1\} \to \R} A^{\bm u}(\widetilde{R})$ 
where $\widetilde{R}(\bEta) = \widetilde{R}_{\bs_0}(\eta_{\bs_0}) \Pi_{\bs\in S_M \setminus \{\bs_0\}} R_\bs(\eta_\bs)$ for all $\bEta = (\eta_\bs)_{\bs\in S_M}$.
\EndFor
$v_{\textsf{new}} = A_M^{\bm u}(R)$
\EndWhile
\end{algorithmic}
\end{algorithm} 
\begin{rmrk}
To compute the constants $\alpha_i$, we can either explicitly implement the partial evaluations of $A^{\bm u}_M(R)$. 
Alternatively, we can treat $A^{\bm u}_M(R)$ as a function in $\R^2\to \R$ depending on the values $R_{\bs_0}(0)$ and $R_{\bs_0}(1)$. We know that this function is a multivariate-polynomial of the form $\alpha_1 R_{\bs_0}(1)^2 + \alpha_2 R_{\bs_0}(0)^2 + \alpha_3 R_{\bs_0}(1)R_j(0) + \alpha_4 R_{\bs_0}(1) + \alpha_5 R_{\bs_0}(0) + \alpha_6$. 
The constants can be computed using  multivariate-polynomial interpolation in six points. This interpolation has the advantages that it is non-intrusive and that the evaluations of $A^{\bm u}_M(R)$ can be performed efficiently using the ideas of Section~\ref{sec:rewritingAu}.
\end{rmrk}

%% file: 5DeterministicSolving.tex
\section{Deterministic resolution of \revision{a} cross-diffusion \revision{system}}
\label{sec:FiniteVolumes}

We describe in this section the numerical scheme we use in Section~\ref{sec:numerical:experiments} for the resolution of \revision{ a cross-diffusion problem, which may be seen as a simplified version of} system~\eqref{eq:cross:diff:system}. \revision{The numerical scheme is baseed on} a cell-centered finite volume method~\cite{Eymard:1999, Eymard:Gallouet:Herbin:2000}, assuming that a numerical approximation of the self-diffusion matrix $\matD_s(\overline{\rho})$ can be computed for any $\overline{\rho}\in [0,1]$.
\revision{The design and numerical analysis of a numerical scheme for the approximation of the original problem\eqref{eq:cross:diff:system} is left for future work.} 

\medskip

\revision{The simplified cross-diffusion system we consider here reads as follows. Let $\Omega \subset \mathbb{R}^d$ be a polyhedric bounded domain of $\mathbb{R}^d$. Local volumic fractions of  $\re$ and $\bl$ particles are given by functions $\rho_{\re}: \Omega \times \mathbb{R}_{+} \to [0,1]$ and $\rho_{\bl}: \Omega \times \mathbb{R}_{+}  \to [0,1]$\revision{, where} $\rho_{\re}(\bx,t)$ (respectively  $\rho_{\bl}(\bx,t)$) denotes the local volumic fraction of $\re$ (respectively $\bl$) particles at point $\bx \in \Omega$ and time $t\geq 0$. We assume here that 
 $\rho_{\re}$ and $\rho_{\bl}$ {are} solutions to the following simplified cross-diffusion system:}
 \revision{
\begin{align}
\label{eq:cross:diff:systemsimp}
    \frac{\partial}{\partial t} 
    \begin{pmatrix} 
    \rho_{\re} \\ \rho_{\bl} 
    \end{pmatrix} 
    = \frac{1}{2} \nabla \cdot \left[
    \begin{pmatrix} 
    \dps {\rm Tr}\left[\frac{\rho_{\bl}}{\rho} \matD_s(\rho) + \dps \frac{\rho_{\re}}{\rho} \matD\right] & \dps {\rm Tr}\left[\frac{\rho_{\re}}{\rho} (\matD-\matD_s(\rho)) \right]\\ 
    \dps{\rm Tr}\left[ \frac{\rho_{\bl}}{\rho} (\matD-\matD_s(\rho)) \right]
    & \dps {\rm Tr}\left[\frac{\rho_{\re}}{\rho} \matD_s(\rho) + \dps \frac{\rho_{\bl}}{\rho} \matD\right] 
    \end{pmatrix} 
    \begin{pmatrix} \nabla \rho_{\re} \\ \nabla \rho_{\bl}
    \end{pmatrix}\right], \quad {\mbox{ in }\Omega \times \mathbb{R}_+^*}, 
\end{align}
{where $\rho:= \rho_\re + \rho_\bl$.} 
We complement system~\eqref{eq:cross:diff:system} by the initial conditions $\rho_{\re}(0,\bx) = \rho_{\re}^0(\bx),\ \rho_{\bl}(0,\bx) = \rho_{\bl}^0(\bx)$, where 
{$\rho_{\re}^0, \rho_\bl^0 \in L^\infty(\Omega)$ satisfying} $0 \leq \rho_{\re}^0+\rho_{\bl}^0 \leq 1$ almost everywhere in $\Omega$}. 
\revision{We also assume here that system (\ref{eq:cross:diff:systemsimp}) is complemented with Neumann boundary conditions.}

\subsection{Notation} 
Let $T_f>0$ denote some final time. 
For the time discretization, we introduce a division of the interval $\left[0,\Tf \right]$ into subintervals $\Ip := \left[t_{p-1},t_{p}\right]$, $1 \leq p \leq \Nt$ for some $\Nt\in \mathbb{N}^*$ such that $0=t_0 < t_1 < \cdots < t_{\Nt}=\Tf$. 
The time steps are denoted by $\Delta t_p=t_{p}-t_{p-1}$, $p=1,\cdots,\Nt$. 
For a given sequence of real numbers $(v^p)_{p\in \mathbb{N}}$, we define the approximation of the first-order time derivative thanks to the backward Euler scheme as follows:
\begin{equation*}
\partial_t v^p := \frac{v^p - v^{p-1}}{\Delta t_p} \quad \forall \,\, 1 \leq p \leq \Nt.
\end{equation*}

\begin{figure}[ht]
    \centering
    \includegraphics[width = 0.32 \textwidth]{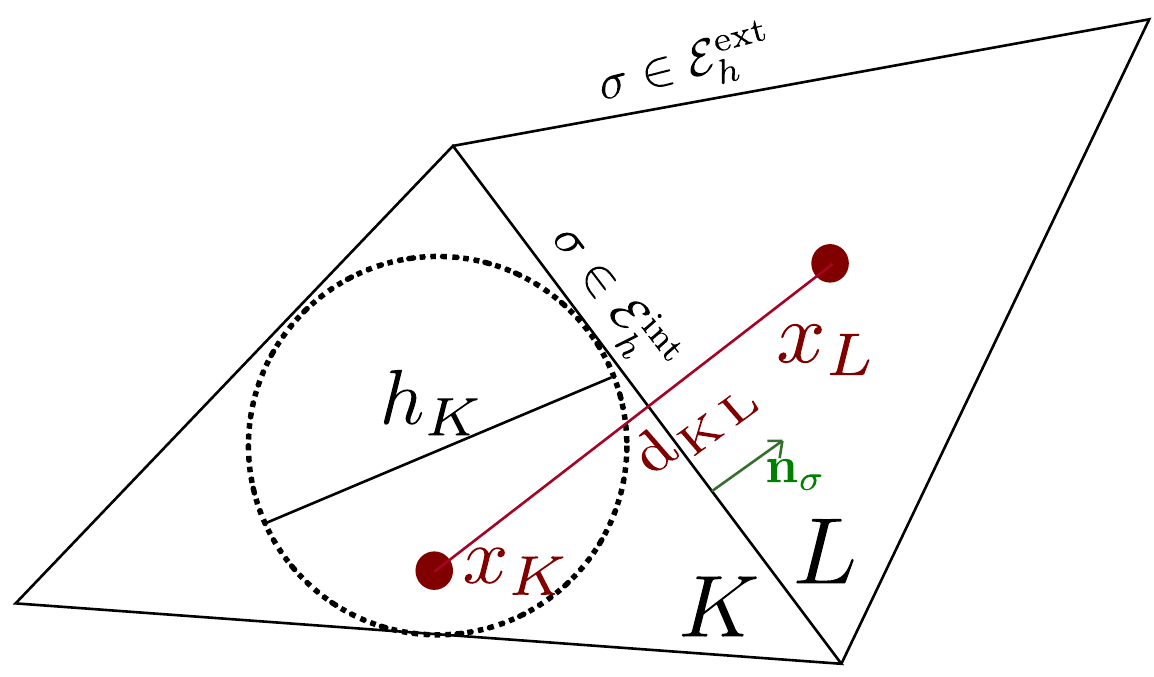}
    \caption{\revision{Illustration of the notation for a mesh with two elements.}}
    \label{ref:fig:mesh:notation}
\end{figure}

For the space discretization, we consider a conforming simplicial mesh $\Th$ of the domain $\Omega$, i.e. $\Th$ is a set of elements $K$ verifying $\dps \bigcup_{K \in \Th} \overline{K}=\overline{\Omega}$,
where the intersection of the closure of two elements of
$\Th$ is either an empty set, a vertex, or a $l$-dimensional face, $0 \leq l \leq d - 1$.
Denote by $\hK$ the diameter of the generic element $K \in \Th$ and  $h \:= \max_{K\in\Th} \hK$.
We denote by $\Eh$ the set of mesh faces.
Boundary faces are collected in the set $\Ehext=\left\{\sigma \in \Eh; \sigma \subset \partial \Omega \right\}$ and internal faces are collected in the set $\Ehint = \Eh \backslash \Ehext$.
To each face $\sigma \in \Eh$, we associate a unit normal vector ${\bm n}_{\sigma}$; for $\sigma \in \Ehint$, $\sigma = K \cap L$, ${\bm n}_{\sigma}$ points from $K$ towards $L$ and for $\sigma \in \Ehext$ it coincides with the outward unit normal vector ${\bm n}_{\Omega}$ of $\Omega$.
We also denote by $\Ne$ the number of elements in the mesh $\Th$. Furthermore, the notation $\bn_{K,\sigma}$ stands for the outward unit normal vector to the element $K$ on $\sigma$.
We also assume that the family $\Th$ is superadmissible in the sense that for
all cells $K \in \Th$ there exists a point $\bx_K \in K$ (the cell
center) and for all edges $\sigma \in \Eh$ there exists a point $\bx_{\sigma} \in \partial K$
(the edge center) such that, for all edges $\sigma \in \EK$, the line
segment joining $\bx_K$ with $\bx_{\sigma}$ is orthogonal to $\sigma$ (see~\cite{Eymard:Gallouet:Herbin:2000}).
For an interior edge $\sigma \in \Ehint$ shared by two elements $K$ and $L$ (denoted in the sequel by $\sigma = K \cap L$), we define the
distance between these elements $d_{KL} := \dist(\bx_K, \bx_L)$. \revision{Figure~\ref{ref:fig:mesh:notation} provides a schematic illustration.}

\subsection{The cell-centered finite volume method}

In the context of the cell-centered finite volume method, the unknowns of the model are discretized using one value per cell: $ \forall 1 \leq p \leq \Nt$ we let
\begin{equation*}
    \bU^p := (\bU_K^p)_{K \in \Th} \in \mathbb{R}^{2 \Ne}, \quad \text{with} \quad \bU_K^p := (\rho_{\re, K}^p, \rho_{\bl, K}^p  ) \in \mathbb{R}^2
\end{equation*}
where $\rho_{\re, K}^p$ and $\rho_{\bl, K}^p$ are respectively the discrete elementwise unknowns approximating the values of $\rho_{\re}(t_p)$ and $\rho_{\bl}(t_p)$ in the element $K \in \Th$.
More precisely,
$\rho_{\re, K}^p \approx \dps \frac{1}{|K|} \int_{K} \rho_{\re}(t_p,x) \,\mathrm{dx}$ and $\rho_{\bl, K}^p \approx \dps \frac{1}{|K|} \int_{K} \rho_{\bl}(t_p,x) \,\mathrm{dx}$.
Integrating~\eqref{eq:cross:diff:system} over the element $K \in \Th$, using next the Green formula and the Neumann boundary conditions, the cell-centered finite volume scheme we solve is the following : for a given fixed value of $\bU^{p-1}$, find $\bU^p \in \mathbb{R}^{2 \Ne}$ such that
\begin{equation}
\label{eq:FV:scheme}
\begin{split}
       & H_{1,K}(\bU^p) = |K|\partial_t \rho_{\re,K}^p - \frac{1}{2} \sum_{\sigma \in \EK} \F_{1,K,\sigma}(\bU^p)=0 \quad \forall K \in \Th,\\
       & H_{2,K}(\bU^p) = |K|\partial_t \rho_{\bl,K}^p - \frac{1}{2} \sum_{\sigma \in \EK} \F_{2,K,\sigma}(\bU^p)=0 \quad \forall K \in \Th.
\end{split}
\end{equation}

The numerical fluxes $\F_{1,K,\sigma}(\bU^p)$ and $\F_{2,K,\sigma}(\bU^p)$ are respectively approximations of
\begin{equation*}
\begin{split}
    \F_{1,K,\sigma}(\bU^p) & 
    \approx \left( \matS^{11}(\bU^p) \nab \rho_{\re}^p +  
\dps \matS^{12}(\bU^p) \nab \rho_{\bl}^p\right) \cdot \bn_{K,\sigma}
\\
    \F_{2,K,\sigma}(\bU^p) & 
    \approx
    \left(\matS^{21}(\bU^p) \nab \rho_{\re}^p +  
\dps \matS^{22}(\bU^p) \nab \rho_{\bl}^p \right) \cdot \bn_{K,\sigma},
\end{split}
\end{equation*}
where $\matS^{11}(\bU^p) := \left(\mathcal{S}^{11}_K(\bU^p)\right)_{K\in \mathcal T_h}$, $\matS^{12}(\bU^p) := \left(\mathcal{S}^{12}_K(\bU^p)\right)_{K\in \mathcal T_h}$, $\matS^{21}(\bU^p) := \left(\mathcal{S}^{21}_K(\bU^p\right))_{K\in \mathcal T_h}$, $\matS^{22}(\bU^p) := \left(\mathcal{S}^{22}_K(\bU^p)\right)_{K\in \mathcal T_h}$ are collections of \revision{real numbers} defined  on each cell of the mesh as follows. For all $K\in \mathcal T_h$,  
\revision{
\begin{equation*}
\begin{split}
    \mathcal{S}^{11}_{K}(\bU^p) &:= \dps {\rm Tr}\left[\frac{\rho_{\bl,K}^p}{\rho_K^p} \matD_s(\rho_K^p) + \dps \frac{\rho_{\re,K}^p}{\rho_K^p} \matD\right], \\
  \mathcal{S}^{12}_{K}(\bU^p) &:= \dps {\rm Tr}\left[\frac{\rho_{\re,K}^p}{\rho_K^p} \left(\matD-\matD_s(\rho_K^p)\right) \right] \\
    \mathcal{S}^{21}_{K}(\bU^p)&:= \dps {\rm Tr}\left[\frac{\rho_{\bl,K}^p}{\rho_K^p} \left(\matD-\matD_s(\rho_K^p)\right) \right] \\
    \mathcal{S}^{22}_{K}(\bU^p) &:= \dps {\rm Tr}\left[\frac{\rho_{\re,K}^p}{\rho_K^p} \matD_s(\rho_K^p) + \dps \frac{\rho_{\bl,K}^p}{\rho_K^p} \matD \right],
    \end{split}
\end{equation*}
}
where $\rho_K^p:= \rho_{\re,K}^p + \rho_{\bl,K}^p$. Now we are interested in computing an approximation of $\left( \matS^{ij}(\bU^p) \nab \rho_{\sharp}^p \right) \cdot \bn_{K,\sigma}$ for all $1\leq i,j \leq 2$ and $\sharp \in \left\{\re,\bl \right\}$. 

For any edge $\sigma \in \Ehint$, denoting by $K,L\in \mathcal T_h$ the two cells sharing the edge $\sigma$, we employ the following harmonic averaging formula for all $1\leq i,j \leq 2$,
\revision{
\begin{equation*}
\left( \matS^{ij}(\bU^p) \nab \rho_{\sharp}^p \right) \cdot \bn_{K,\sigma} \approx 
|\sigma| \frac{\mathcal{S}^{ij}_{K}(\bU^p) \mathcal{S}^{ij}_{L}(\bU^p)}{d_{K,\sigma} \mathcal{S}_{L}^{ij}(\bU^p) + d_{L,\sigma} \mathcal{S}^{ij}_{K}(\bU^p)} \left( \rho_{\sharp,L}^p - \rho_{\sharp,K}^p \right).
\end{equation*}
}
We refer to~\cite{Eymard:Herbin:2004,Eymard:Herbin:2007,Agelas:Eymard:Herbin:2009} for more details. This averaging formula enables then to obtain the expression of the fluxes $ \F_{1,K,\sigma}(\bU^p) $ and $ \F_{2,K,\sigma}(\bU^p) $.

At each time step $p$, we thus have to solve the nonlinear problem: find $\bU^p\in \mathbb{R}^{2\Ne}$
\begin{equation*}
    \mathcal{H}^p(\bU^p) = 0 
\end{equation*}
where $\mathcal{H}^p(\bU^p):= \left( H_{1,K}(\bU^p), H_{2,K}(\bU^p)\right)_{K\in \mathcal T_h}\in \mathbb{R}^{2 \Ne}$, where 
$H_{1,K}(\bU^p)$ and $H_{2,K}(\bU^p)$ are defined by~\eqref{eq:FV:scheme}.

\subsection{\revision{Resolution of the nonlinear problem by a Newton method}}
\label{ref:sec:Newton:methods}
In this section, we present a Newton procedure for computing the approximate solution $\bU^p$ at each time step $p$.
For $1 \leq p \leq \Nt$ and $\bU^{p,0} \in \mathbb{R}^{2 \Ne}$ fixed (typically $\bU^{p,0} = \bU^{p-1}$), the Newton algorithm generates a
sequence $(\bU^{p,k})_{k \geq 1}$, with $\bU^{p,k} \in \mathbb{R}^{2 \Ne}$ given by the system of linear algebraic equations:
\begin{equation}
\label{eq:Newton:scheme}
    \matA^{p,k-1} \bU^{p,k} = \bB^{p,k-1}.
\end{equation}
Here, the Jacobian matrix $\matA^{p,k-1} \in \mathbb{R}^{2 \Ne, 2 \Ne}$ and the right-hand side vector $\bB^{p,k-1} \in \mathbb{R}^{2 \Ne}$ are defined by
\begin{equation*}
    \matA^{p,k-1}:=\matJ_{\mathcal{H}^p}(\bU^{p,k-1}) \quad \text{and} \quad  \bB^{p,k-1}:= \matJ_{\mathcal{H}^p}(\bU^{p,k-1}) \bU^{p,k-1} - \mathcal{H}^p(\bU^{p,k-1}).
\end{equation*}
Note that $\matJ_{\mathcal{H}^p}(\bU^{p,k-1})$ is the Jacobian matrix of the function $\mathcal{H}^p$ at point $\bU^{p,k-1}$.
A classical stopping criterion for system~\eqref{eq:Newton:scheme} is for instance
\begin{equation}
\label{Newton:stopping:criterion}
    \left\|\mathcal{H}^p(\bU^{p,k})\right\|_2 /     \left\|\mathcal{H}^p(\bU^{p,0})\right\|_2 < \epslin
\end{equation}
where $\epslin > 0$ is small enough.
We refer to the book~\cite{Kelley:1995} for large descriptions on linearization techniques.

%% file: 6NumericalExperiments.tex
\section{Numerical Experiments}\label{sec:numerical:experiments}
In the following numerical experiments~\footnote{{The code to reproduce all experiments is available on \url{https://github.com/cstroessner/SelfDiffusionCoefficent.git}}}, we consider the periodic lattice with {$M=1$} and $d=2$, which implies $N=8$. The possible displacement vectors are $\bv_1=(1,0)$, $\bv_2=(-1,0)$, $\bv_3=(0, 1)$, and $\bv_4=(0,-1)$ so that $K=4$ with associated probability $p_k=1/4$ for $k=1,2,3,4$. Throughout this section, we set the tolerance in Algorithm~\ref{alg:ALS} to $\eps = 10^{-5}$. {Note that using the low-rank approximation described in Section~\ref{sec:low:rank} for larger values of $M$ is currently work in progress. }

\begin{rmrk} \label{rem:Interpolation}
Equation~\eqref{eq:approx:Ds} can only be evaluated for certain values {$\overline{\rho} = \frac{\ell}{N}$ for some $0\leq \ell \leq N$}. To evaluate $\mathbb{D}^M_s(\overline{\rho})$ for general $\overline{\rho} \in [0,1]$, we interpolate the individual entries of $\mathbb{D}^M_s(\overline{\rho})$ based on the nodes, which are admissible for Equation~\eqref{eq:approx:Ds}, using splines.
\end{rmrk}

\subsection{{Approximation} of the self-diffusion coefficient}
\subsubsection{Low-rank approximation error}\label{sec:low:rank:error:experiments}
We are interested in studying the effect of minimizing~\eqref{eq:def:Au} with respect to functions $\Psi \in H_M$ compared to functions in $R \in \mathcal T_M^1$. {We denote by } $\Psi_{\textsf{LSQ}}^{\bm u}\in \mathbb{R}^{2^N}$ the solution to (\ref{eq:linlsq}) computed numerically by solving the linear least squares problem in Section~\ref{sec:relaxing} and the low-rank solution $R_{\textsf{ALS}}^{\bm u}$ {obtained} using Algorithm~\ref{alg:ALS}. 
Note that $R_{\textsf{ALS}}^{\bm u}$ may vary depending on the initial choice of $R^0_\bs$ for $\bs\in S_M$. Thus, we run Algorithm~\ref{alg:ALS} $100$ times with different initializations for $R_\bs(0)$ and $R_\bs(1)$ obtained by sampling from the uniform distribution on $[0,1]$. 
For $\bm u = (1,0)$ we obtain $A^{\bm u}_M(\Psi^{\bm u}_{\textsf{LSQ}}) \approx \revision{53.594}$ and $\text{mean}(A^{\bm u}_M(R^{\bm u}_{\textsf{ALS}})) = \revision{53.759}$ with $\min(A^{\bm u}_M(R^{\bm u}_{\textsf{ALS}})) = \revision{53.758}$ and $\max(A^{\bm u}_M(R^{\bm u}_{\textsf{ALS}})) = \revision{53.766}$. 
These results show that evaluating $A_M^{\bm u}(R^{\bm u}_{\textsf{ALS}})$ differs from $A_M^{\bm u}(\Psi^{\bm u}_{\textsf{LSQ}})$ by less than $0.\revision{3}\%$. 
\revision{At this point, we want to emphasize that the goal of Algorithm~\ref{alg:ALS} is to compute a pure tensor product minimizer of $A_M^{\bm u}$, which can potentially be different from a pure tensor product approximation of $\Psi^{\bm u}_{\textsf{LSQ}}$.}

For the computation of the self-diffusion coefficient, we need to evaluate $A_{M,\ell}^{\bm u}$ instead of $A_M^{\bm u}$. We again consider the approximations $\Psi^{\bm u}_{\textsf{LSQ}},R^{\bm u}_{\textsf{ALS}}$, for which we obtain the values displayed in Table~\ref{tab:Alu:experiment}. We observe that the values obtained from $R^{\bm u}_{\textsf{ALS}}$ are again very close to the ones obtained from $\Psi^{\bm u}_{\textsf{LSQ}}$.

\begin{table}[ht]
    \centering
    {\tablinesep=5pt\tabcolsep=5pt
\begin{tabular}{ll|lllllllll}
 &  & \multicolumn{1}{l|}{$\ell=0$} & \multicolumn{1}{l|}{$\ell=1$} & \multicolumn{1}{l|}{$\ell=2$} & \multicolumn{1}{l|}{$\ell=3$} & \multicolumn{1}{l|}{$\ell=4$} & \multicolumn{1}{l|}{$\ell=5$} & \multicolumn{1}{l|}{$\ell=6$} & \multicolumn{1}{l|}{$\ell=7$} & \multicolumn{1}{l|}{$\ell=8$}  \\ \hline
\multicolumn{1}{l|}{\multirow{2}{*}{$\bm u = \dps \binom{1}{0}$}} & $A_{\revision{M,}\ell}^{\bm u}(\Psi^{\bm u}_{\textsf{LSQ}})$  &      
 \revision{0.5000}  &  0.4196  &  0.3430  &  0.2708 &    0.2035 &    0.1421 &    0.0873 &    0.0398    &     0
    \\ \cline{2-2}
\multicolumn{1}{l|}{}                           & $A_{\revision{M,}\ell}^{\bm u}(R^{\bm u}_{\textsf{ALS}})$ & 
  0.5000   & 0.4197 &   0.3433  &  0.2714  &  0.2044 &    0.1430 &    0.0878 &    0.0398    &     0   \\ \cline{1-2}
\multicolumn{1}{l|}{\multirow{2}{*}{$\bm u =        \dps \binom{0}{1}$}} & $A_{\revision{M,}\ell}^{\bm u}(\Psi^{\bm u}_{\textsf{LSQ}})$   &
 0.5000   & 0.4197 &   0.3430 &   0.2708&    0.2035 &    0.1421 &   0.0873 &    0.0398 &         0
\\ \cline{2-2}
\multicolumn{1}{l|}{}                           & $A_{\revision{M,}\ell}^{\bm u}(R^{\bm u}_{\textsf{ALS}})$ &       
0.5000   & 0.4197 &   0.3433    & 0.2714    & 0.2044&    0.1430 &    0.0878 &   0.0398 &         0
       \\ \cline{1-2}
\multicolumn{1}{l|}{\multirow{2}{*}{$\bm u = \dps \binom{1}{1}$}} & $A_{\revision{M,}\ell}^{\bm u}(\Psi^{\bm u}_{\textsf{LSQ}})$  &       1.0000  &  0.8393 &    0.6860 &    0.5416 &    0.4070 &    0.2843 &   0.1747  &   0.0795  &       0  \\ \cline{2-2}
\multicolumn{1}{l|}{}                           & $A_{\revision{M,}\ell}^{\bm u}(R^{\bm u}_{\textsf{ALS}})$ &       
1.0000   & 0.8393  &   0.6866  &   0.5428  &  0.4089 &    0.2860 &    0.1756 &    0.0796 &         \revision{0} \\ \cline{1-2}
\end{tabular}
  \label{tab:Alu:experiment}
}
    \caption{Evaluation of $A_{\revision{M},\ell}^{\bm u}$ for $\Psi^{\bm u}_{\textsf{LSQ}}$ and $R^{\bm u}_{\textsf{ALS}}$ for different values of $\bm u$ and $\ell$.}
\end{table}

\subsubsection{Computing the self-diffusion coefficient}
\begin{figure}[ht]
    \centering
    \includegraphics[width = 0.4\textwidth]{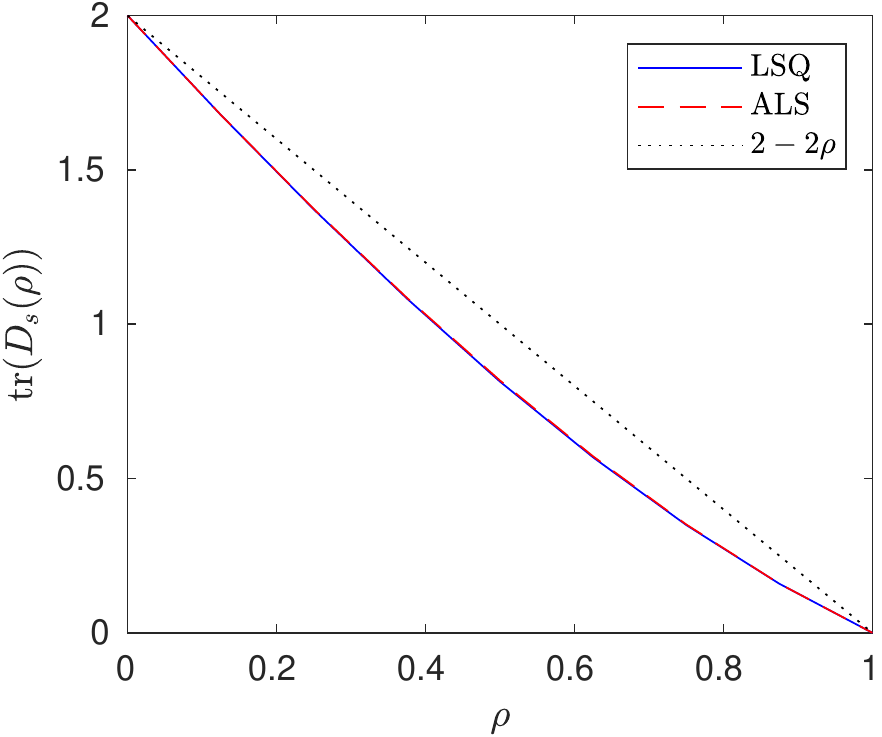}
    \caption{Plot of the trace of $\mathbb D_s(\rho)$ computed based on $\Psi_{\textsf{LSQ}}^{\bm u}$ and $R_{\textsf{ALS}}^{\bm u}$ and plot of the function \revision{$2-2\rho$}. }
    \label{fig:trace:8}
\end{figure}

For $0\leq \ell \leq N$ we evaluate equation~\eqref{eq:approx:Ds} for $\bm u \in \{(1,0),(0,1),(1,1)\}$ to recover the entries of the symmetric matrix $\mathbb D_s\left(\frac{\ell}{N}\right)$. 
To evaluate the self-diffusion coefficient $\mathbb{D}_s(\overline{\rho})$ for different values of $\overline{\rho}\in[0,1]$, we use an interpolation technique as explained in Remark~\ref{rem:Interpolation}. \revision{We observe that the symmetry in the studied jumping scheme leads to off-diagonal entries, which are numerically zero. Additionally, we find that $(\mathbb{D}_s(\overline{\rho}))_{11} = (\mathbb{D}_s(\overline{\rho}))_{22}$. In Figure~\ref{fig:trace:8}, we plot ${\rm Tr}\left(\mathbb D_s(\overline{\rho})\right)$,  where ${\rm Tr}$ \revision{denotes} the trace operator.  Our approximation satisfies the known properties $\mathbb D_s(0) =  \mathbb I$  and $\mathbb D_s(1) =  0\cdot \mathbb I$, where $\mathbb I$ denotes the identity matrix.
}
Note that both trace plots almost coincide.

\subsection{Resolution of the cross-diffusion system}
In the following, we solve the PDE system~\eqref{eq:cross:diff:system} based on the self-diffusion coefficient obtained from $R^{\bm u}_{\textsf{ALS}}$. 
The system is defined on the unit square domain $\Omega :=(0, 1) \times (0, 1)$. 
We consider a uniform spatial mesh ($\Ne = 312$ elements) and we use a constant time step $\Delta t_p = \Delta t = 10^{-3}$ seconds $\forall 1 \leq p \leq \Nt$.
The final time of simulation is $\Tf = 10$ seconds.
We employ at each time step $1 \leq p \leq \Nt$ a Newton solver as described in Section~\ref{ref:sec:Newton:methods} with $\epslin = 10^{-8}$.
The initial values are defined by
\begin{equation*}
    \rho_{\re}^{0}(x,y) := 0.25 + 0.25 \cos(\pi x) \cos(\pi y) \quad \text{and} \quad \rho_{\bl}^{0}(x,y) := 0.5 -  0.5 \cos(\pi x) \cos(\pi y).
\end{equation*}
\begin{figure}[ht]
    \centering
    \includegraphics[width = 0.32 \textwidth]{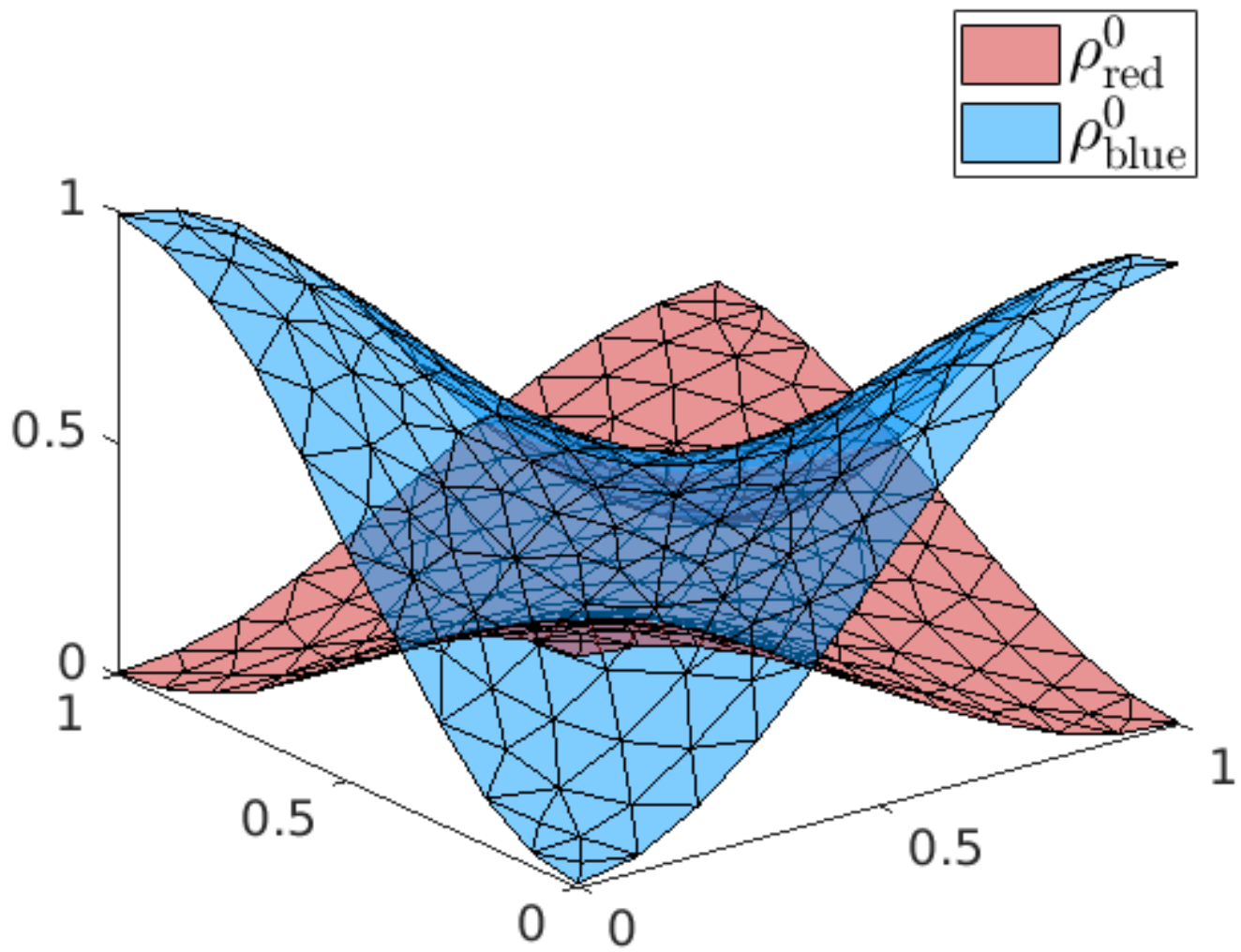}
        \includegraphics[width = 0.32 \textwidth]{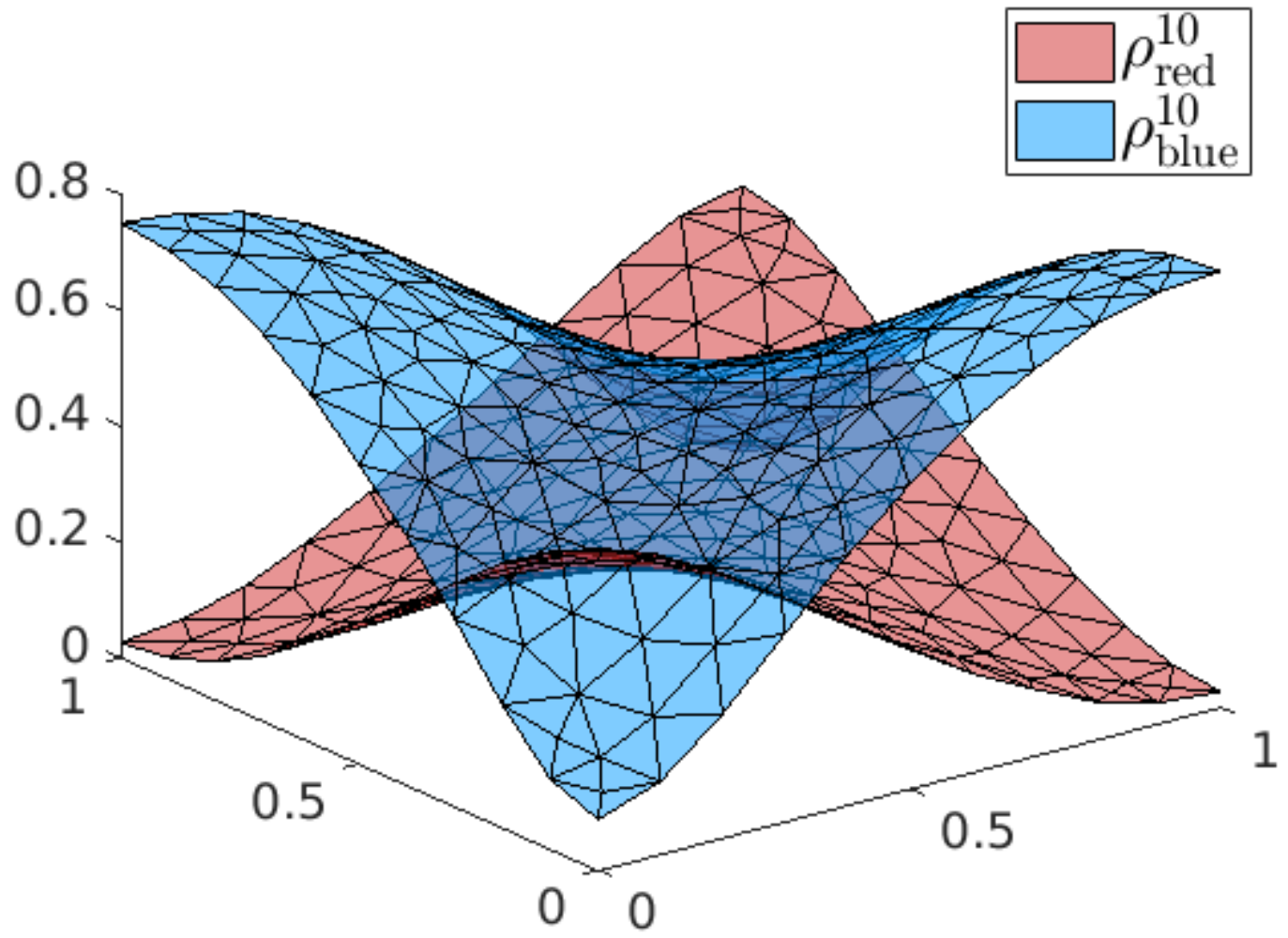}
    \includegraphics[width = 0.32 \textwidth]{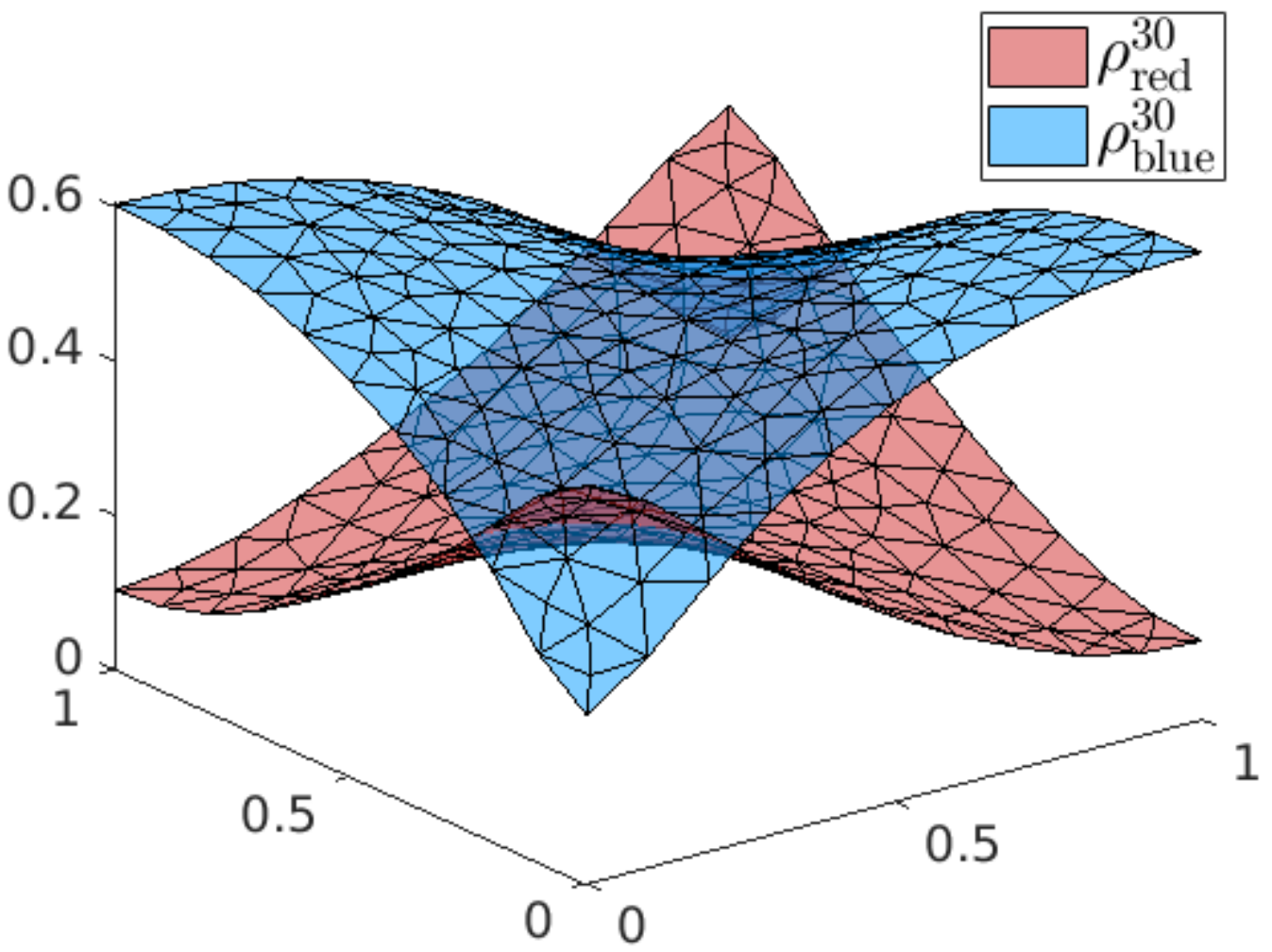}
    \includegraphics[width = 0.32 \textwidth]{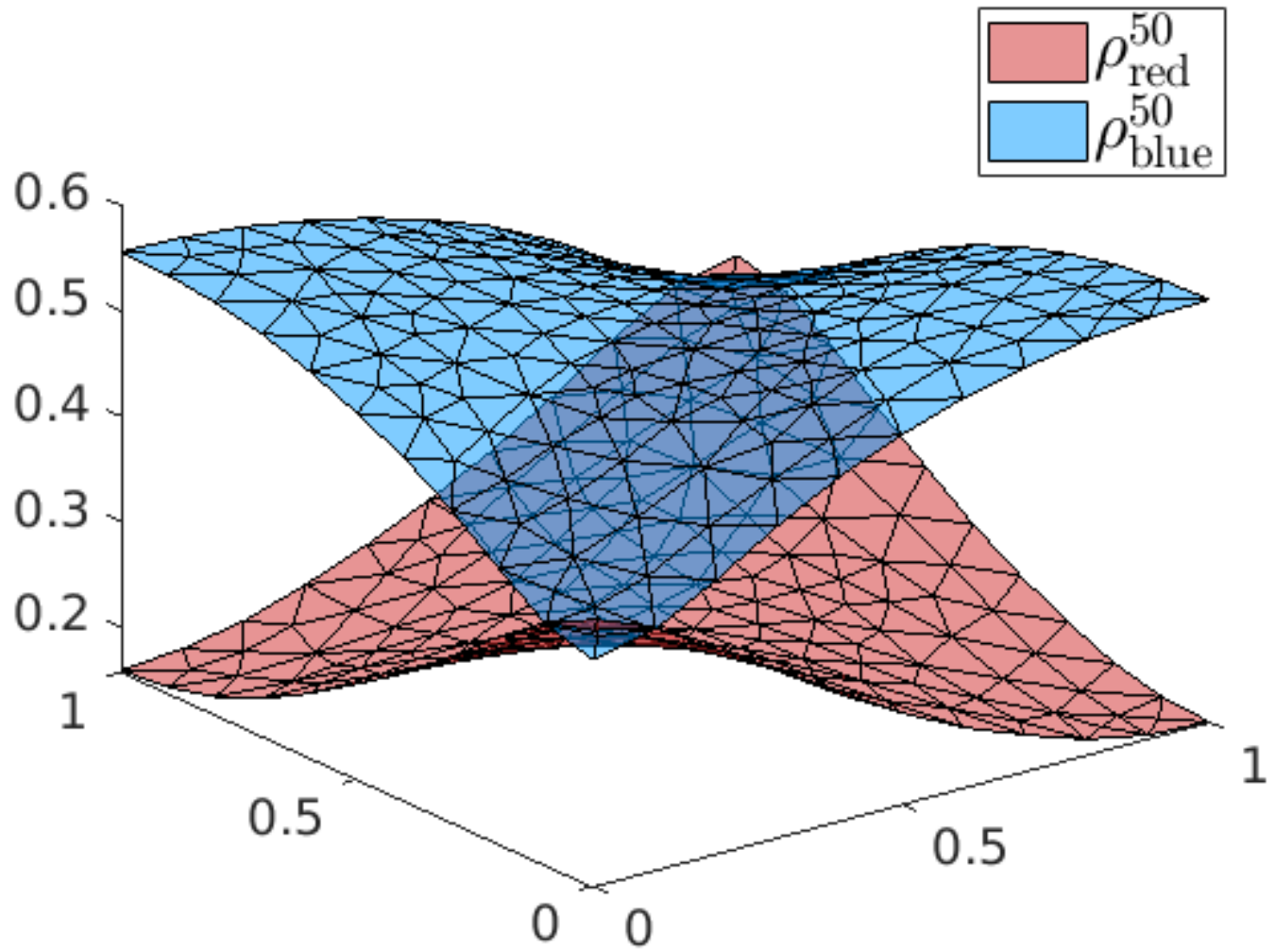}
     \includegraphics[width = 0.32 \textwidth]{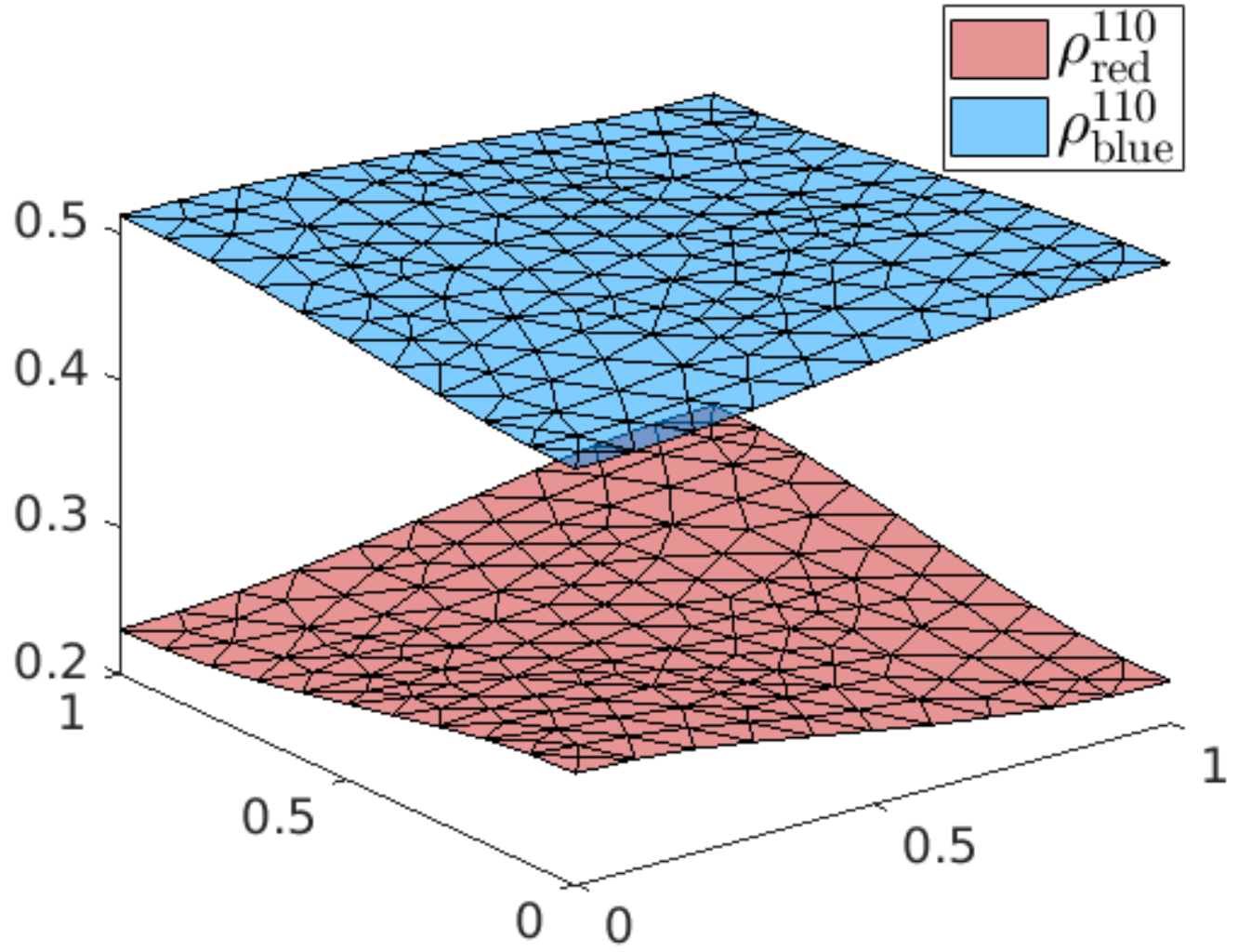}
    \includegraphics[width = 0.32 \textwidth]{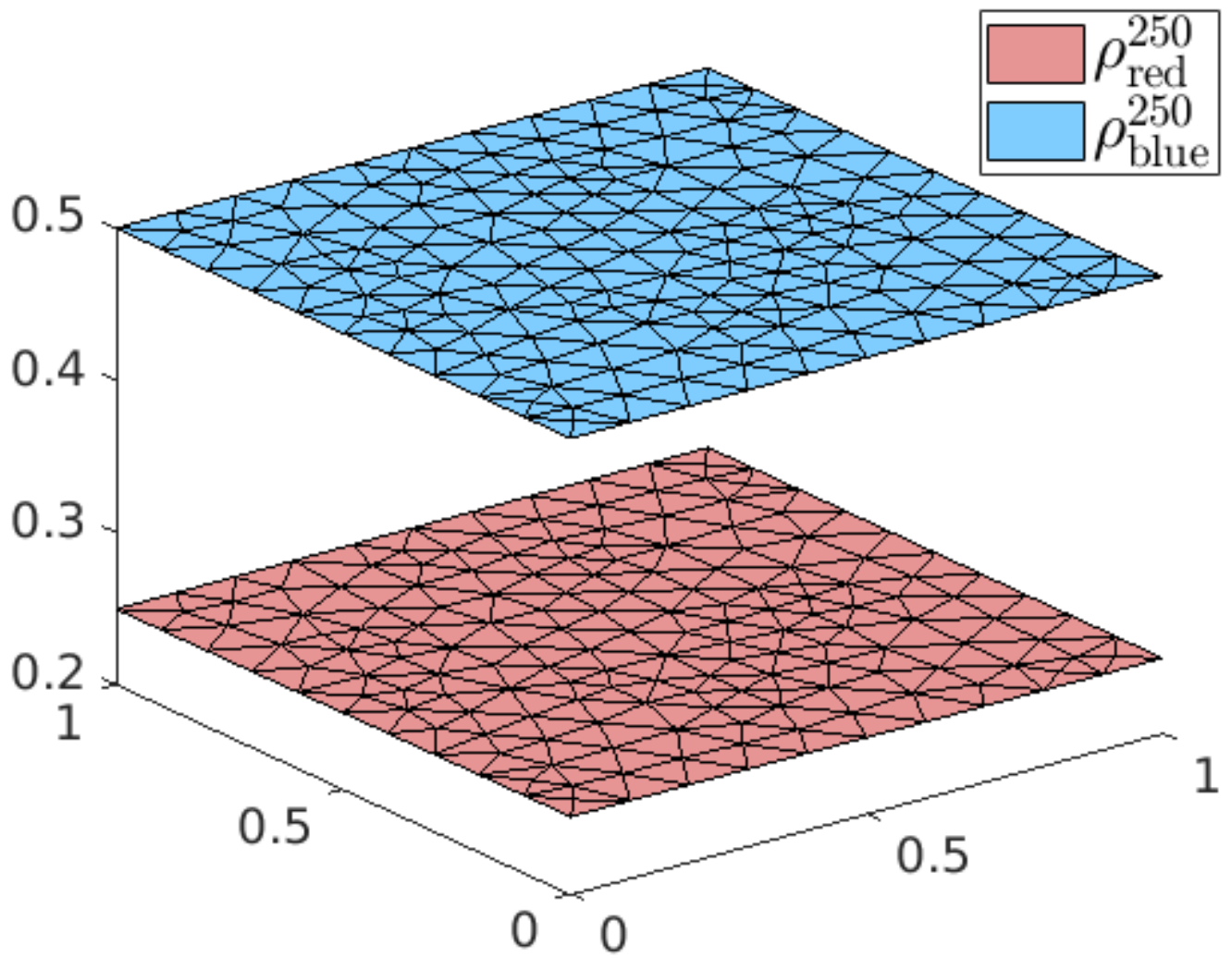}
    \caption{Solutions $\rho_{\re}$ and $\rho_{\bl}$ for several time steps. Top left $p=0$, top middle $p=10$, top right $p=30$, bottom left $p=50$, bottom middle $p=110$, bottom right $p=250$.}
    \label{ref:fig:solutions:several:times}
\end{figure}
Figure~\ref{ref:fig:solutions:several:times} displays the shape of the numerical solutions $\rho_{\re}$ and $\rho_{\bl}$ for several time steps. 
We observe that the local volumic fractions evolve over time to reach constant profiles (around $0.25$ for $\rho_{\re}$ and $0.5$ for $\rho_{\bl}$) in the long time limit. \revision{This indicates, that our scheme is stable in the sense that the densities lie in their physical ranges. In every time step, only $2$ or $3$ Newton iterations are required to reach the stopping criterion~\eqref{Newton:stopping:criterion}.}

%% file: 7Conclusion.tex
\section{Conclusion}
In this work, we derived a novel way of computing the self-diffusion {matrix of the tagged particle process} via low-rank solutions of a high-dimensional optimization problem. {The obtained approximation can then in turn be used in order to compute the solution of cross-diffusion systems that arise as the hydrodynamic limits of multi-species symmetric exclusion systems like the one introduced in~\cite{Quastel92}, } using a finite volume scheme. 
Our numerical results have clearly demonstrated that the computed low-rank solutions led to accurate approximations of the self-diffusion coefficient, {at least for small finite lattice sizes}. 
\revision{The e}xtension to \revision{more sophisticated low-rank approximation formats and} larger lattices, for which solving the minimization problem directly is intractable,  as well as the derivation of a convergent entropy diminishing finite volume scheme \revision{for hydrodynamic limits of multi-species exclusion processes} will be explored in a future work.